\documentclass[11pt]{article}
\usepackage{cite}
\usepackage{mathrsfs}
\usepackage{amsfonts}
\usepackage{amsmath}
\usepackage{amsfonts,amssymb}
\usepackage{dsfont}
\usepackage{curves}
\usepackage{mathrsfs}
\usepackage{pifont}
\usepackage{amssymb}
\allowdisplaybreaks

\numberwithin{equation}{section}

\date{}

\def\BigRoman{\uppercase\expandafter{\romannumeral\number\count 255 }}
\def\Romannumeral{\afterassignment\BigRoman\count255=}

\begin{document}
\title{Spanning $k$-trees and distance spectral radius in graphs
%\thanks{Supported by }
}
\author{\small  Sizhong Zhou\footnote{Corresponding
author. E-mail address: zsz\_cumt@163.com (S. Zhou), wujiancheng@just.edu.cn (J. Wu)}, Jiancheng Wu\\
\small  School of Science, Jiangsu University of Science and Technology,\\
\small  Zhenjiang, Jiangsu 212100, China\\
}

\maketitle
\begin{abstract}
\noindent Let $k\geq2$ be an integer. A tree $T$ is called a $k$-tree if $d_T(v)\leq k$ for each $v\in V(T)$, that is, the maximum degree of a $k$-tree
is at most $k$. A $k$-tree $T$ is a spanning $k$-tree if $T$ is a spanning subgraph of a connected graph $G$. Let $\lambda_1(D(G))$ denote the distance
spectral radius in $G$, where $D(G)$ denotes the distance matrix of $G$. In this paper, we verify an upper bound for $\lambda_1(D(G))$ in a connected
graph $G$ to guarantee the existence of a spanning $k$-tree in $G$.
\\
\begin{flushleft}

(2020) Mathematics Subject Classification: 05C50, 05C05, 05C70
\end{flushleft}
\end{abstract}

\section{Introduction}

Throughout this paper, we deal only with finite and undirected graphs which admit neither multiple edges nor loops. Let $G$ be a graph with
vertex set $V(G)=\{v_1,v_2,\ldots,v_n\}$ and edge set $E(G)$. The degree of a vertex $v$ in $G$, written $d_G(v)$, is the number of vertices
adjacent to $v$. We denote by $c(G)$ the number of connected components in $G$. For a given subset $S\subseteq V(G)$, we use $G[S]$ and $G-S$
to denote the subgraphs of $G$ induced by $S$ and $V(G)\setminus S$, respectively. The number of elements in the set $S$ is denoted by $|S|$.
The path and complete graph of order $n$ are denoted by $P_n$ and $K_n$, respectively.

The distance between $v_i$ and $v_j$ in $G$, denoted by $d_G(v_i,v_j)$ (or $d_{ij}$), is the length of a shortest path from $v_i$ to $v_j$. The
distance matrix of $G$, denoted by $D(G)$, is the $n$-by-$n$ real symmetric matrix whose $(i,j)$-entry is $d_G(v_i,v_j)$ (or $d_{ij}$). The
distance eigenvalues in $G$ are the eigenvalues of its distance matrix $D(G)$. Let $\lambda_1(D(G))\geq\lambda_2(D(G))\geq\cdots\geq\lambda_n(D(G))$
be its distance eigenvalues in nonincreasing order. In view of the Perron-Frobenius theorem, $\lambda_1(D(G))$ is always positive (unless $G$
is trivial), and we call $\lambda_1(D(G))$ the distance spectral radius in $G$.

For two graphs $G_1$ and $G_2$, we call $G_1$ and $G_2$ isomorphic, and write $G_1\cong G_2$, if there exists a bijection $\varphi: V(G_1)\rightarrow
V(G_2)$ with $uv\in E(G_1)\Leftrightarrow \varphi(u)\varphi(v)\in E(G_2)$ for all $u,v\in V(G_1)$. We denote the disjoint union of $G_1$ and $G_2$
by $G_1\cup G_2$, which is the graph with vertex set $V(G_1\cup G_2)=V(G_1)\cup V(G_2)$ and edge set $E(G_1\cup G_2)=E(G_1)\cup E(G_2)$. In particular,
if $G_1\cong G_2$, then write $2G_1=G_1\cup G_2$ for short. We denote the join of $G_1$ and $G_2$ by $G_1\vee G_2$, which is the graph with vertex
set $V(G_1\vee G_2)=V(G_1)\cup V(G_2)$ and edge set $E(G_1\vee G_2)=E(G_1)\cup E(G_2)\cup\{vu:u\in V(G_1), v\in V(G_2)\}$.

For two integers $a$ and $b$ such that $0\leq a\leq b$, a spanning subgraph $F$ of $G$ is called an $[a,b]$-factor if $a\leq d_F(v)\leq b$ for
each $v\in V(G)$. If $a=1$, then an $[a,b]$-factor is a $[1,b]$-factor. If $b=1$, then a $[1,b]$-factor is a 1-factor or a perfect matching.
Let $k\geq2$ be an integer. A tree $T$ is called a $k$-tree if $d_T(v)\leq k$ for each $v\in V(T)$, that is, the maximum degree of a $k$-tree
is at most $k$. A $k$-tree $T$ is a spanning $k$-tree if $T$ is a spanning subgraph of a connected graph $G$. It is obvious that $G$ admits a
spanning $k$-tree if and only if $G$ admits a connected $[1,k]$-factor.

Brouwer and Haemers \cite{BH} described some sufficient conditions, by means of the adjacency eigenvalues and Laplacian eigenvalues, for a
graph to have a perfect matching. Up to now, much attention has been paid to this topic, and we refer the reader to \cite{CGH,LZ,O,ZLp}. Lots
of researchers posed some sufficient conditions to guarantee the existence of $[1,2]$-factors in graphs by using various parameters, such as
the degree conditions \cite{AEKKM,Zd}, the binding number \cite{L}, the independence number \cite{Wp,ZSL2,WZi}, the isolated toughness
\cite{GCW,Zs} and others \cite{ZZS,ZSL1,ZWB}. Some results on $[a,b]$-factors or fractional
$[a,b]$-factors of graphs were derived by Kouider \cite{K}, Kouider and Lonc \cite{KL}, Wang and Zhang \cite{WZs}, Wu \cite{Wa}, Zhou et al.
\cite{ZL,ZPX,ZPX1,Za1,ZZs,Zr,ZZL}. Ding, Johnson and Seymour \cite{DJS}, Gargano, Hammar, Hell, Stacho and Vaccaro \cite{GHHSV} studied the
existence of a spanning tree with a given parameter in a connected graph. Kyaw \cite{Kyaw} showed a degree
and neighborhood condition for a connected graph to admit a spanning $k$-tree. Neumann-Lara and Rivera-Campo \cite{NR} presented an independence
number condition for the existence of a spanning $k$-tree in a connected graph.

Very recently, Gu and Liu \cite{GL} discovered a sufficient condition for a spanning $k$-tree in a connected graph by using Laplacian eigenvalues.
Fan, Goryainov, Huang and Lin \cite{FGHL} obtained two sufficient conditions for a connected graph to possess a spanning $k$-tree by using
adjacency spectral radius and signless Laplacian spectral radius, respectively. Along this line, we extends the existing body of knowledge on
spanning $k$-trees by focusing on the distance spectral radius, a less traditional measure compared to the usual adjacency or Laplacian spectral
radius. This novel approach can enrich the spectral graph theory and its applications in studying the structural properties of graphs. The study
of spanning $k$-trees is significant in various applications, such as network design, where ensuring certain connectivity properties is crucial.
The paper's results can provide new tools for analyzing and designing such networks.

In what follows, we put forward a distance spectral radius condition which guarantees the existence of a spanning $k$-tree in a connected graph.

\medskip

\noindent{\textbf{Theorem 1.1.}} Let $k\geq4$ be an integer, and let $G$ be a connected graph of order $n\geq k+2$.

(\romannumeral1) For $k=4$ and $n\geq12$, or $k\geq5$, if $\lambda_1(D(G))\leq\lambda_1(D(G^{*}))$, then $G$ admits a spanning $k$-tree unless
$G\cong G^{*}$, where $G^{*}\cong K_1\vee(K_{n-k-1}\cup kK_1)$.

(\romannumeral2) For $k=4$ and $6\leq n\leq11$, if $\lambda_1(D(G))\leq\lambda_1(D(G^{\sharp}))$, then $G$ admits a spanning $4$-tree unless
$G\cong G^{\sharp}$, where $G^{\sharp}\cong K_{\frac{n-3}{3}}\vee(\frac{2n+3}{3}K_1)$.

\section{Preliminary lemmas}

In this section, we put forward several necessary preliminary lemmas, which are very useful to prove our main result. A fundamental property
of the distance spectral radius of a graph and its spanning subgraph, which is a corollary of the Perron-Frobenius theorem, was derived in the book
\cite{G}.

\medskip

\noindent{\textbf{Lemma 2.1}} (Godsil \cite{G}). Let $e$ be an edge of a graph $G$ such that $G-e$ is connected. Then
$$
\lambda_1(D(G))<\lambda_1(D(G-e)).
$$

We now explain the concepts of quotient matrices and equitable partitions.

\noindent{\textbf{Definition 2.2}} (Brouwer and Haemers \cite{BH1}). Let $M$ be a real matrix of order $n$ described in the following block form

\begin{align*}
M=\left(
  \begin{array}{ccc}
    M_{11} & \cdots & M_{1r}\\
    \vdots & \ddots & \vdots\\
    M_{r1} & \cdots & M_{rr}\\
  \end{array}
\right),
\end{align*}
where the blocks $M_{ij}$ are $n_i\times n_j$ matrices for any $1\leq i,j\leq r$ and $n=n_1+n_2+\cdots+n_r$. For $1\leq i,j\leq r$, let $b_{ij}$
denote the average row sum of $M_{ij}$, that is, $b_{ij}$ is the sum of all entries in $M_{ij}$ divided by the number of rows. Then $B(M)=(b_{ij})$
(simply by $B$) is called a quotient matrix of $M$. If for every pair $i,j$, $M_{ij}$ admits constant row sum, then the partition is called
equitable.

\noindent{\textbf{Lemma 2.3}} (You, Yang, So and Xi \cite{YYSX}). Let $M$ be a real matrix with an equitable partition, and let $B$ be the
corresponding quotient matrix. Then the eigenvalues of $B$ are also eigenvalues of $M$. Furthermore, if $M$ is a nonnegative matrix, then
$\rho_1(B)=\rho_1(M)$, where $\rho_1(B)$ and $\rho_1(M)$ denote the largest eigenvalues of the matrices $B$ and $M$, respectively.

\medskip

\noindent{\textbf{Lemma 2.4}} (Zhang and Lin \cite{ZLp}). Let $n,t,s$ and $n_i$ ($1\leq i\leq t$) be positive integers with
$n_1\geq n_2\geq\cdots\geq n_t\geq1$ and $n_1+n_2+\cdots+n_t=n-s$. Then
$$
\lambda_1(D(K_s\vee(K_{n_1}\cup K_{n_2}\cup\cdots\cup K_{n_t})))\geq\lambda_1(D(K_s\vee(K_{n-s-(t-1)}\cup(t-1)K_1))),
$$
with equality if and only if $K_s\vee(K_{n_1}\cup K_{n_2}\cup\cdots\cup K_{n_t})\cong K_s\vee(K_{n-s-(t-1)}\cup(t-1)K_1)$.

\medskip

In 1989, Win \cite{W} verified the following result, which shows a sufficient condition for the existence of a spanning $k$-tree in a connected
graph. Ellingham and Zha \cite{EZ} posed a shorter proof.

\medskip

\noindent{\textbf{Lemma 2.5}} (Win \cite{W}, Ellingham and Zha \cite{EZ}). Let $k\geq3$ be an integer. If a connected graph $G$ satisfies
$$
c(G-S)\leq(k-2)|S|+2,
$$
for any vertex subset $S$ of $G$, then $G$ admits a spanning $k$-tree.

\medskip

Let $W(G)=\sum_{i<j}d_{ij}$ be the Wiener index of a connected graph $G$ of order $n$. Note that
$\lambda_1(D(G))=\max\limits_{X\in\mathbb{R}^{n}}\frac{X^{T}D(G)X}{X^{T}X}$. Then we possess
$$
\lambda_1(D(G))=\max_{X\in\mathbb{R}^{n}}\frac{X^{T}D(G)X}{X^{T}X}\geq\frac{\mathbf{1}^{T}D(G)\mathbf{1}}{\mathbf{1}^{T}\mathbf{1}}\geq\frac{2W(G)}{n},
$$
where $\mathbf{1}=(1,1,\ldots,1)^{T}$.

\section{The proof of Theorem 1.1}

In what follows, we verify Theorem 1.1.

\noindent{\it Proof of Theorem 1.1.} Assume to the contrary that a connected graph $G$ has no spanning $k$-tree with the minimum distance spectral
radius. In terms of Lemma 2.5, we derive
$$
c(G-S)\geq(k-2)|S|+3,
$$
for some nonempty subset $S$ of $V(G)$. Together with Lemma 2.1 and the choice of $G$, we have that the induced subgraph $G[S]$ and all connected
components in $G-S$ are complete graphs. Furthermore, we deduce $G=G[S]\vee(G-S)$. Let $|S|=s$ and $c(G-S)=t$. Clearly,
$G\cong K_s\vee(K_{n_1}\cup K_{n_2}\cup\cdots\cup K_{n_t})$
for some positive integers $n_1\geq n_2\geq\cdots\geq n_t$ with $n_1+n_2+\cdots+n_t=n-s$. Let $n_2=n_3=\cdots=n_t=1$. Then we obtain a graph
$G'\cong K_s\vee(K_{n-s-(t-1)}\cup(t-1)K_1)$. By Lemma 2.4, we see $\lambda_1(D(G'))\leq\lambda_1(D(G))$ with equality if and only if
$G\cong G'$. Thus, we infer $G\cong G'$ (otherwise $\lambda_1(D(G'))<\lambda_1(D(G))$ which contradicts that $G$ has the minimum distance spectral
radius). Recall that $t\geq(k-2)s+3$. Let $\widetilde{G}\cong K_s\vee(K_{n-(k-1)s-2}\cup((k-2)s+2)K_1)$. In what follows, we compare the distance
spectral radius of $G'$ and $\widetilde{G}$.

\noindent{\bf Claim 1.} $\lambda_1(D(\widetilde{G}))\leq\lambda_1(D(G'))$ with equality if and only if $G'\cong\widetilde{G}$.

\noindent{\it Proof.} If $t=(k-2)s+3$, then $G'\cong\widetilde{G}$. Next, we consider $t\geq(k-2)s+4$. We denote the vertex set of $\widetilde{G}$
by $V(\widetilde{G})=V(K_s)\cup V(K_{n-(k-1)s-2})\cup V(((k-2)s+2)K_1)$. Assume that $X$ is the Perron vector of $D(\widetilde{G})$, and let $X(v)$
denote the entry of $X$ corresponding to the vertex $v\in V(\widetilde{G})$. In terms of symmetry, it is obvious that all vertices of $K_s$ (resp.
$K_{n-(k-1)s-2}$ and $((k-2)s+2)K_1$) admit the same entries in $X$. Thus, we may let $X(u)=x_1$ for any vertex $u$ in $((k-2)s+2)K_1$, $X(v)=x_2$
for any vertex $v$ in $K_{n-(k-1)s-2}$, and $X(w)=x_3$ for any vertex $w$ in $K_s$. It follows from $t\geq(k-2)s+4$ and $n\geq s+t\geq (k-1)s+4$
that
\begin{align*}
&\lambda_1(D(G'))-\lambda_1(D(\widetilde{G}))\geq X^{t}(D(G')-D(\widetilde{G}))X\\
=&(t-(k-2)s-4)(t-(k-2)s-3)x_2^{2}+(t-(k-2)s-3)(n-s-(t-1))x_2^{2}\\
&+(t-(k-2)s-3)(n-s-(t-1))x_2^{2}\\
>&0,
\end{align*}
which implies that $\lambda_1(D(\widetilde{G}))<\lambda_1(D(G'))$ if $t\geq(k-2)s+4$. This completes the proof of Claim 1. \hfill $\Box$

Recall that $G\cong G'$. Combining this with Claim 1, we get $G\cong\widetilde{G}$ (otherwise $\lambda_1(D(\widetilde{G}))<\lambda_1(D(G'))=\lambda_1(D(G))$
which contradicts that $G$ admits the minimum distance spectral radius). Let $G^{*}\cong K_1\vee(K_{n-k-1}\cup kK_1)$. In the remainder of this section,
we shall clarify that, in most cases, the graph $G^{*}$ admits the minimum distance spectral radius and contains no spanning $k$-tree.

\noindent{\bf Claim 2.} If $n\geq(k-1)s+4$, then $\lambda_1(D(G^{*}))\leq\lambda_1(D(\widetilde{G}))$ with equality if and only if $\widetilde{G}\cong G^{*}$.

\noindent{\it Proof.} Recall that $\widetilde{G}\cong K_s\vee(K_{n-(k-1)s-2}\cup((k-2)s+2)K_1)$. Then the quotient matrix of distance matrix
$D(\widetilde{G})$ by virtue of the partition $\{V(K_s), V(K_{n-(k-1)s-2}), V(((k-2)s+2)K_1)\}$ can be expressed as
\begin{align*}
\left(
  \begin{array}{ccc}
    s-1 & n-(k-1)s-2 & (k-2)s+2\\
    s & n-(k-1)s-3 & 2(k-2)s+4\\
    s & 2n-2(k-1)s-4 & 2(k-2)s+2\\
  \end{array}
\right).
\end{align*}
The corresponding characteristic polynomial is
\begin{align*}
f(x)=&x^{3}-(n+(k-2)s-2)x^{2}\\
&-(2(k-2)sn+7n-(k-2)(2k-1)s^{2}-(7k-8)s-11)x\\
&+(k-2)s^{2}n-2(k-3)sn-6n-(k-2)(k-1)s^{3}\\
&+(2k^{2}-9k+8)s^{2}+2(4k-7)s+10.
\end{align*}

If $s=1$, then $\widetilde{G}\cong G^{*}$ and the polynomial $f(x)$ becomes $g(x)=x^{3}-(n+k-4)x^{2}-((2k+3)n-2k^{2}-2k-5)x
-(k+2)n+k^{2}+2k+2$. In view of Lemma 2.3, $\lambda_1(D(\widetilde{G}))$ is the largest root of $f(x)=0$ and $\lambda_1(D(G^{*}))=\theta(n)$ (simply
$\theta$) is the largest root of $g(x)=0$. In the following, we let $s\geq2$ so that $n\geq(k-1)s+4\geq2k+2$. Note that $g(\theta)=0$. By plugging
the value $\theta$ into $x$ of $f(x)-g(x)$, we obtain
\begin{align*}
f(\theta)=&f(\theta)-g(\theta)\\
=&(s-1)(-(k-2)\theta^{2}+(-2(k-2)n+(k-2)(2k-1)s+2k^{2}+2k-6)\theta\\
&+(k-2)sn-(k-4)n-(k-2)(k-1)s^{2}+(k^{2}-6k+6)s+k^{2}+2k-8).
\end{align*}

Set $p(\theta)=-(k-2)\theta^{2}+(-2(k-2)n+(k-2)(2k-1)s+2k^{2}+2k-6)\theta+(k-2)sn-(k-4)n-(k-2)(k-1)s^{2}+(k^{2}-6k+6)s+k^{2}+2k-8$. Note that
\begin{align*}
\theta=&\lambda_1(D(G^{*}))\geq\frac{2W(G^{*})}{n}\\
=&\frac{n^{2}+(2k-1)n-k^{2}-3k}{n}\\
=&\frac{n^{2}+(k-1)n+kn-k^{2}-3k}{n}\\
\geq&\frac{n^{2}+(k-1)n+k^{2}-k}{n}\\
>&n+k-1,
\end{align*}
and so $\frac{-2(k-2)n+(k-2)(2k-1)s+2k^{2}+2k-6}{2(k-2)}=-n+\frac{(2k-1)s}{2}+\frac{k^{2}+k-3}{k-2}<n+k-1<\theta$. Thus, we deduce
$$
p(\theta)<p(n+k-1).
$$

Let $h(n)=p(n+k-1)=-3(k-2)n^{2}+(2k(k-2)s-2k^{2}+13k-10)n-(k-2)(k-1)s^{2}+(2k^{3}-6k^{2}+k+4)s+k^{3}+5k^{2}-11k$. Note that $n\geq(k-1)s+4$ and
$\frac{2k(k-2)s-2k^{2}+13k-10}{6(k-2)}<(k-1)s+4\leq n$. Thus, we infer
\begin{align*}
h(n)\leq&h((k-1)s+4)\\
=&-(k-1)(k-2)^{2}s^{2}-(7k^{2}-34k+34)s+k^{3}-3k^{2}-7k+56.
\end{align*}

Let $q(s)=h((k-1)s+4)=-(k-1)(k-2)^{2}s^{2}-(7k^{2}-34k+34)s+k^{3}-3k^{2}-7k+56$. Then
\begin{align*}
q'(s)=&-2(k-1)(k-2)^{2}s-(7k^{2}-34k+34)\\
\leq&-4(k-1)(k-2)^{2}-(7k^{2}-34k+34)\\
=&-4k^{3}+13k^{2}+2k-18\\
=&-k((4k-1)(k-3)-5)-18\\
<&0,
\end{align*}
by $s\geq2$ and $k\geq4$. Thus, we have
\begin{align*}
q(s)\leq&q(2)=-4(k-1)(k-2)^{2}-2(7k^{2}-34k+34)+k^{3}-3k^{2}-7k+56\\
=&-3k^{3}+3k^{2}+29k+4\\
=&-k((3k-12)(k+3)+7)+4\\
\leq&-7k+4\\
<&0,
\end{align*}
by $k\geq4$. So when $n\geq(k-1)s+4$, we infer
\begin{align*}
f(\theta)=&(s-1)p(\theta)<(s-1)p(n+k-1)\\
=&(s-1)h(n)\leq(s-1)h((k-1)s+4)\\
=&(s-1)q(s)<0,
\end{align*}
which implies that $\lambda_1(D(G^{*}))=\theta<\lambda_1(D(\widetilde{G}))$. We finish the proof of Claim 2. \hfill $\Box$

\noindent{\bf Claim 3.} If $n=(k-1)s+3$ and $k\geq5$, then $\lambda_1(D(G^{*}))\leq\lambda_1(D(\widetilde{G}))$ with equality if and only if
$\widetilde{G}\cong G^{*}$.

\noindent{\it Proof.} Recall that $\widetilde{G}\cong K_s\vee(K_{n-(k-1)s-2}\cup((k-2)s+2)K_1)$ and $G^{*}\cong K_1\vee(K_{n-k-1}\cup kK_1)$. If
$s=1$, then $\widetilde{G}\cong G^{*}$. In what follows, we consider $s\geq2$. Then we infer
\begin{align*}
h(n)=&-3(k-2)n^{2}+(2k(k-2)s-2k^{2}+13k-10)n\\
&-(k-2)(k-1)s^{2}+(2k^{3}-6k^{2}+k+4)s+k^{3}+5k^{2}-11k\\
=&-3(k-2)((k-1)s+3)^{2}+(2k(k-2)s-2k^{2}+13k-10)((k-1)s+3)\\
&-(k-2)(k-1)s^{2}+(2k^{3}-6k^{2}+k+4)s+k^{3}+5k^{2}-11k\\
=&-(k-1)(k-2)^{2}s^{2}-(3k^{2}-20k+22)s+k^{3}-k^{2}+k+24.
\end{align*}

Let $r(s)=-(k-1)(k-2)^{2}s^{2}-(3k^{2}-20k+22)s+k^{3}-k^{2}+k+24$. Then
\begin{align*}
r'(s)=&-2(k-1)(k-2)^{2}s-(3k^{2}-20k+22)\\
\leq&-4(k-1)(k-2)^{2}-(3k^{2}-20k+22)\\
=&-4k^{3}+17k^{2}-12k-6\\
=&-k((4k-13)(k-1)-1)-6\\
<&0,
\end{align*}
by $s\geq2$ and $k\geq5$. Thus, we derive
\begin{align*}
h(n)=&r(s)\leq r(2)=-4(k-1)(k-2)^{2}-2(3k^{2}-20k+22)+k^{3}-k^{2}+k+24\\
=&-3k^{3}+13k^{2}+9k-4\\
=&-k((3k+2)(k-5)+1)-4\\
<&0,
\end{align*}
by $k\geq5$. Hence, we get
\begin{align*}
f(\theta)=&(s-1)p(\theta)<(s-1)p(n+k-1)\\
=&(s-1)h(n)<0,
\end{align*}
which implies that $\lambda_1(D(G^{*}))=\theta<\lambda_1(D(\widetilde{G}))$. Claim 3 is verified. \hfill $\Box$

Next, we consider $n=(k-1)s+3$ and $k=4$. Observe that $\widetilde{G}\cong G^{\sharp}$, where
$G^{\sharp}\cong K_{\frac{n-3}{3}}\vee\left(\frac{2n+3}{3}K_1\right)$. Then the quotient matrix of distance matrix $D(G^{\sharp})$ in terms of
the partition $\left\{V\left(K_{\frac{n-3}{3}}\right),V\left(\frac{2n+3}{3}K_1\right)\right\}$ can be expressed as
\begin{align*}
\left(
  \begin{array}{ccc}
    \frac{n-6}{3} & \frac{2n+3}{3}\\
    \frac{n-3}{3} & \frac{4n}{3} \\
  \end{array}
\right).
\end{align*}
The corresponding characteristic polynomial is
$$
\varphi(x)=x^{2}-\frac{5n-6}{3}x+\frac{2n^{2}-21n+9}{9}.
$$
By virtue of Lemma 2.3, the largest root of $\varphi(x)=0$ equals the distance spectral radius $\lambda_1(D(G^{\sharp}))$ of
$G^{\sharp}\cong K_{\frac{n-3}{3}}\vee\left(\frac{2n+3}{3}K_1\right)$, namely,
$$
\rho(G^{\sharp})=\frac{5n-6+\sqrt{17n^{2}+24n}}{6}.
$$

\noindent{\bf Claim 4.} If $n=(k-1)s+3=3s+3$ and $n\geq12$, then $\lambda_1(D(G^{*}))<\lambda_1(D(G^{\sharp}))$.

\noindent{\it Proof.} Let $n=(k-1)s+3=3s+3$ and $n\geq13$, we can compute $h(n)=-\frac{4}{3}n^{2}+\frac{34}{3}n+54<0$. Hence, we infer
$f(\theta)=(s-1)p(\theta)<(s-1)p(n+k-1)=(s-1)h(n)<0$, which implies that $\lambda_1(D(G^{*}))=\theta<\lambda_1(D(\widetilde{G}))=\lambda_1(D(G^{\sharp}))$.

Note that $k=4$ and $\lambda_1(D(G^{*}))=\theta(n)$ is the largest root of $g(x)=x^{3}-nx^{2}-(11n-45)x-6n+26=0$. For $s=3$ and $n=(k-1)s+3=12$, we have
$g(x)=x^{3}-12x^{2}-87x-46$ and $\lambda_1(D(G^{*}))=\lambda_1(D(K_3\vee9K_1))=9+2\sqrt{19}$. It is clear that $g(9+2\sqrt{19})=68+32\sqrt{19}>0$ and
$g'(9+2\sqrt{19})=168+60\sqrt{19}>0$, so $\lambda_1(D(G^{\sharp}))>\theta(12)=\lambda_1(D(G^{*}))$. This completes the proof of Claim 4. \hfill $\Box$

\noindent{\bf Claim 5.} If $n=(k-1)s+3=3s+3$ and $6\leq n\leq11$, then $\lambda_1(D(G^{*}))\geq\lambda_1(D(G^{\sharp}))$ with equality if and only if
$G^{\sharp}\cong G^{*}$.

\noindent{\it Proof.} For $s=1$ and $n=k+2=6$, we get $G^{\sharp}\cong G^{*}$. For $s=2$ and $n=2k+1=9$, we derive $\lambda_1(D(G^{\sharp}))=\lambda_1(D(K_2\vee7K_1))=\frac{13+\sqrt{177}}{2}$, $g(x)=x^{3}-9x^{2}-54x-28$, and $g\left(\frac{13+\sqrt{177}}{2}\right)=-20<0$,
so $\lambda_1(D(G^{\sharp}))<\theta(9)=\lambda_1(D(G^{*}))$. Claim 5 is proved. \hfill $\Box$

From the discussion above, we always admit $G\cong G^{\sharp}$ if $k=4$ and $6\leq n\leq11$ and $G\cong G^{*}$ if $k=4$ and $n\geq12$, or $k\geq5$
and $n\geq k+2$, otherwise we always possess $\lambda_1(D(G^{\sharp}))<\lambda_1(D(G))$ if $k=4$ and $6\leq n\leq11$ and
$\lambda_1(D(G^{*}))<\lambda_1(D(G))$ if $k=4$ and $n\geq12$, or $k\geq5$ and $n\geq k+2$, which contradicts the hypothesis of this theorem. This
completes the proof of Theorem 1.1. \hfill $\Box$

\section*{Data availability statement}

My manuscript has no associated data.

\section*{Declaration of competing interest}

The authors declare that they have no conflicts of interest to this work.

\medskip

\section*{Acknowledgments}

The authors would like to show their great gratitude to anonymous referees for the valuable suggestions which largely improve the quality
and the readability of this paper. Project ZR2023MA078 supported by Shandong Provincial Natural Science Foundation.


\begin{thebibliography}{9999}

\bibitem {BH} A. Brouwer, W. Haemers, Eigenvalues and perfect matchings, Linear Algebra and its Applications 395 (2005) 155--162.

\bibitem {CGH} S. Cioab\v{a}, D. Gregory, W. Haemers, Matchings in regular graphs from eigenvalue, Journal of Combinatorial Theory, Series B 99 (2009) 287--297.

\bibitem {LZ} S. Li, M. Zhang, On the signless Laplacian index of cacti with a given number of pendant vertices, Linear Algebra and its
Applications 436 (2012) 4400--4411.

\bibitem {O} S. O, Spectral radius and matchings in graphs, Linear Algebra and its Applications 614 (2021) 316--324.

\bibitem {ZLp} Y. Zhang, H. Lin, Perfect matching and distance spectral radius in graphs and bipartite graphs, Discrete Applied Mathematics 304 (2021) 315--322.

\bibitem {AEKKM} K. Ando, Y. Egawa, A. Kaneko, K. Kawarabayashi, H. Matsuda, Path factors in claw-free graphs, Discrete Mathematics 243 (2002) 195--200.

\bibitem {Zd} S. Zhou, Degree conditions and path factors with inclusion or exclusion properties, Bulletin Mathematique de la Societe des Sciences
Mathematiques de Roumanie 66(1) (2023) 3--14.

\bibitem {L} H. Liu, Binding number for path-factor uniform graphs, Proceedings of the Romanian Academy, Series A: Mathematics, Physics, Technical
Sciences, Information Science 23(1) (2022) 25--32.

\bibitem{ZSL2} S. Zhou, Z. Sun, H. Liu, Some sufficient conditions for path-factor uniform graphs, Aequationes Mathematicae 97(3) (2023) 489--500.

\bibitem {WZi} S. Wang, W. Zhang, Independence number, minimum degree and path-factors in graphs, Proceedings of the Romanian Academy, Series A:
Mathematics, Physics, Technical Sciences, Information Science 23(3) (2022) 229--234.

\bibitem {Wp} J. Wu, Path-factor critical covered graphs and path-factor uniform graphs, RAIRO-Operations Research 56(6) (2022) 4317--4325.

\bibitem{Zs} S. Zhou, Some results on path-factor critical avoidable graphs, Discussiones Mathematicae Graph Theory 43(1) (2023) 233--244.

\bibitem {GCW} W. Gao, Y. Chen, Y. Wang, Network vulnerability parameter and results on two surfaces, International
Journal of Intelligent Systems 36 (2021) 4392--4414.

\bibitem {ZZS} S. Zhou, Y. Zhang, Z. Sun, The $A_{\alpha}$-spectral radius for path-factors in graphs, Discrete Mathematics 347(5) (2024) 113940.

\bibitem {ZSL1} S. Zhou, Z. Sun, H. Liu, Distance signless Laplacian spectral radius for the existence of path-factors in graphs, Aequationes
Mathematicae 98(3) (2024) 727--737.

\bibitem {ZWB} S. Zhou, J. Wu, Q. Bian, On path-factor critical deleted (or covered) graphs, Aequationes Mathematicae 96(4) (2022) 795--802.

\bibitem {K} M. Kouider, Sufficient condition for the existence of an even $[a,b]$-factor in graph, Graphs and Combinatorics 29(4) (2013) 1051--1057.

\bibitem {KL} M. Kouider, Z. Lonc, Stability number and $[a,b]$-factors in graphs, Journal of Graph Theory 46(4) (2004) 254--264.

\bibitem {WZs} S. Wang, W. Zhang, Some results on star-factor deleted graphs, Filomat 38(3) (2024) 1101--1107.

\bibitem {Wa} J. Wu, A sufficient condition for the existence of fractional $(g,f,n)$-critical covered graphs, Filomat 38(6) (2024) 2177--2183.

\bibitem {ZL} S. Zhou, H. Liu, Two sufficient conditions for odd $[1,b]$-factors in graphs, Linear Algebra and its Applications 661 (2023) 149--162.

\bibitem{Za1} S. Zhou, A neighborhood union condition for fractional $(a,b,k)$-critical covered graphs, Discrete Applied Mathematics
323 (2022) 343--348.

\bibitem {Zr} S. Zhou, Remarks on restricted fractional $(g,f)$-factors in graphs, Discrete Applied Mathematics 354 (2024) 271--278.

\bibitem{ZZs} S. Zhou, Y. Zhang, Sufficient conditions for fractional $[a,b]$-deleted graphs, Indian Journal of Pure and Applied Mathematics,
https://doi.org/10.1007/s13226-024-00564-w

\bibitem{ZPX} S. Zhou, Q. Pan, L. Xu, Isolated toughness for fractional $(2,b,k)$-critical covered graphs, Proceedings of the Romanian Academy,
Series A: Mathematics, Physics, Technical Sciences, Information Science 24(1) (2023) 11--18.

\bibitem{ZPX1} S. Zhou, Q. Pan, Y. Xu, A new result on orthogonal factorizations in networks, Filomat 38(20)(2024), in press.

\bibitem {ZZL} S. Zhou, Y. Zhang, H. Liu, Some properties of $(a,b,k)$-critical graphs, Filomat 38(16)(2024), in press.

\bibitem {DJS} G. Ding, T. Johnson, P. Seymour, Spanning trees with many leaves, Journal of Graph Theory 37(4) (2001) 189--197.

\bibitem {GHHSV} L. Gargano, M. Hammar, P. Hell, L. Stacho, U. Vaccaro, Spanning spiders and light-splitting switches, Discrete Mathematics 285(1--3) (2004) 83--95.

\bibitem {Kyaw} A. Kyaw, A sufficient condition for a graph to have a $k$-tree, Graphs and Combinatorics 17 (2001) 113--121.

\bibitem {NR} V. Neumann-Lara, E. Rivera-Campo, Spanning trees with bounded degrees, Combinatorica 11(1) (1991) 55--61.

\bibitem {GL} X. Gu, M. Liu, A tight lower bound on the matching number of graphs via Laplacian eigenvalues, European Journal of Combinatorics 101 (2022) 103468.

\bibitem {FGHL} D. Fan, S. Goryainov, X. Huang, H. Lin, The spanning $k$-trees, perfect matchings and spectral radius of graphs, Linear and Multilinear Algebra,
https://doi.org/10.1080/03081087.2021.1985055

\bibitem {G} C. Godsil, Algebraic Combinatorics, Chapman and Hall Mathematics Series, New York, 1993.

\bibitem {BH1} A. Brouwer, W. Haemers, Spectra of Graphs - Monograph, Springer, 2011.

\bibitem {YYSX} L. You, M. Yang, W. So, W. Xi, On the spectrum of an equitable quotient matrix and its application, Linear Algebra and its Applications
577 (2019) 21--40.

\bibitem {W} S. Win, On a connection between the existence of $k$-trees and the toughness of a graph, Graphs and Combinatorics 5 (1989) 201--205.

\bibitem {EZ} M. Ellingham, X. Zha, Toughness, trees, and walks, Journal of Graph Theory 33(3) (2000) 125--137.


\end{thebibliography}
\end{document}